\documentclass[11pt,twoside]{article} 
\usepackage{latexsym} 
\input{amssym.def} 
\input{amssym} 
\newfont{\fra}{eufm10 scaled 1095} 
\newfont{\Bb}{msbm10 scaled 1095} 
\newfont{\Bbg}{msbm10 scaled 1680} 
 
\newcommand\RR{{\mbox{\Bb R}}}

\newcommand\fg{{\frak{g}}}

\newcommand\fl{{\frak l}}

\newcommand\fk{{\frak k}} 
\newcommand\fp{{\frak p}}
\newcommand\fa{{\frak a}} 
\newcommand\fd{{\frak d}}

\newcommand\cT{{\cal T}} 
\newcommand\cN{{\cal N}}
\newcommand\ph{\varphi}

\newcommand{\fgl}{\mathop{{\frak g \frak l}}} 
\newcommand{\fsl}{\mathop{{\frak s \frak l}}}

\newcommand{\fso}{\mathop{{\frak s \frak o}}} 
\newcommand{\fo}{\mathop{{\frak o}}} 
\newcommand{\fsu}{\mathop{{\frak s \frak u}}} 

\newcommand{\Iso}{\mathop{{\rm Iso}}} 
\newcommand{\GL}{\mathop{{\rm GL}}} 
\newcommand{\grO}{{\rm O}}
\newcommand{\Hom}{\mathop{{\rm Hom}}} 
\newcommand{\Aut}{\mathop{{\rm Aut}}} 
\newcommand{\Id}{{{\rm id}}} 
\newcommand{\ad}{{{\rm ad}}} 
 
\newcommand{\Ad}{\mathop{{\rm Ad}}}

\newcommand{\Ker}{\mathop{{\rm ker}}} 
\newcommand{\im}{\mathop{{\rm im}}}

\renewcommand{\Im}{\mathop{{\rm Im}}}

\newcommand{\Span}{{{\rm span}}} 
 
\newcommand{\proj}{{{\rm pr}}} 
\newcommand\ip{{\langle\cdot \,,\cdot \rangle}} 
\newcommand\lb{{[\cdot\,,\cdot]}} 
 
\newcommand\dd{\fd(\fl,\theta_\fl,\fa)} 

\newcommand\proof{{\sl Proof. }} 
\newcommand{\qed}{\hspace*{\fill}\hbox{$\Box$}\vspace{2ex}} 
\newtheorem{theo}{Theorem}[section] 
\newtheorem{pr}[theo]{Proposition}
\newtheorem{de}[theo]{Definition}
\newtheorem{ex}[theo]{Example}
\newtheorem{re}[theo]{Remark}
\newtheorem{co}[theo]{Corollary}
\newtheorem{lm}[theo]{Lemma}
\begin{document} 
\title{Indefinite extrinsic symmetric spaces I} 
\author{Ines Kath}
\maketitle 
\begin{abstract}
We find a one-to-one correspondence between full extrinsic symmetric spaces in (possibly degenerate) inner product spaces and certain algebraic objects called (weak) extrinsic symmetric triples. In particular, this yields a description of arbitrary extrinsic symmetric spaces in pseudo-Euclidean spaces by corresponding infinitesimal objects.\\[2ex]

\end{abstract}
\section{Introduction} 
In this paper we will study submanifolds of pseudo-Euclidean spaces. All these submanifolds are assumed to be non-degenerate. This shall mean that the restriction of the  pseudo-Euclidean product to any tangent space of the submanifold is non-degenerate. A
non-degenerate submanifold $M$ of a pseudo-Euclidean space $V$ is called \emph{symmetric submanifold} or \emph{extrinsic symmetric space} if it is invariant under the reflection at each of its affine normal spaces. In particular, each extrinsic symmetric space is a usual (abstract) symmetric space. Similar to usual symmetric spaces extrinsic symmetric 
spaces in can be characterised by curvature. They are exactly those connected complete submanifolds whose second fundamental form is parallel~\cite{S}. 

Two extrinsic symmetric spaces $M\hookrightarrow V$ and $M'\hookrightarrow V'$ are called isometric, if there is an isometry $V\rightarrow V'$ that maps $M$ to $M'$. Here we will be mainly interested in the classification of isometry classes of extrinsic symmetric spaces.

A series of examples of extrinsic symmetric spaces can be obtained in the following way, see e.g.~\cite{F3,EH,Nai}. Start with a pseudo-Riemannian symmetric space $N=G/K$, where $G$ is the transvection group of $N$ and consider the associated decomposition $\fg=\fk\oplus\fp$ of the Lie algebra $\fg$ of $G$. Suppose that there is an element $\xi$ in $\fp$ such that $\ad(\xi)^3=-\ad(\xi)$ on $\fg$. Then the orbit $\Ad(K)\xi\subset\fp$ is an extrinsic symmetric space.  

Conversely, Ferus could prove that  for any extrinsic symmetric space $M$ in a Euclidean space $V$ the following holds. There is an orthogonal splitting $V=V_1\oplus V_2$, $M=M_1\times M_2$ such that $M_1\hookrightarrow V_1$ is a compact extrinsic symmetric space, which lies in a round sphere and $M_2\hookrightarrow V_2$ is an affine subspace.  We may assume that $M_1\hookrightarrow V_1$ full, i.e. that it is not contained in a proper affine subspace of $V_1$. A nice construction due 
to Ferus shows that $M_1\hookrightarrow V_1$ arises in the above described manner as an orbit of the isotropy representation of a semisimple symmetric space. A classification of extrinsic symmetric spaces in Euclidean spaces now follows from the classification of compact symmetric spaces with a suitable element $\xi$ due to Kobayashi and 
Nagano \cite{KoNa, Na}. 

Let us explain the idea of Ferus' construction. Let $M\hookrightarrow \RR^n$ be a full compact extrinsic symmetric space. Assume that $M$ is contained in the sphere. Fix a point $\xi\in M$. Let $K$ be the connected component of the group generated by all reflections at (affine) normal spaces of $M$. This group is called transvection group of the extrinsic symmetric space $M\hookrightarrow \RR^n$. Since $M$ is full this group acts effectively on $M$. Hence it is isomorphic  to the transvection group of the (abstract) symmetric space $M$. In particular, the metric on $M$ is defined by an invariant positive definite scalar product $\ip_\fk$ on the Lie algebra $\fk$ of $K$. Let us  consider the vector space $\fg=\fk\oplus \RR^n$. Then $\ip_\fg=\ip_\fk\oplus \ip_{{\Bbb R}^n}$ is a scalar product on $\fg$.  Ferus showed that one can extend the Lie bracket of $\fk$ and the natural action of $\fk$ on $\RR^n$ to a Lie bracket of $\fg$ such that $\ip_\fg$ becomes invariant. In this way $\fg$ becomes the Lie algebra of the transvection group of a `big' symmetric space and $M$ is the orbit through $\xi$ of the isotropy representation of this space. The inner derivation $\ad(\xi)$ on $\fg$ defined by $\xi$ has the property $\ad(\xi)^3=-\ad(\xi)$. 

Now let us assume that the ambient  space is pseudo-Euclidean. This situation is much more  involved. The easiest case is that of a full and normal extrinsic symmetric space, where normal means that the intersection of all (affine) normal spaces is not empty. As already noticed by Naitoh \cite{Nai} Ferus' construction can be generalised to this case. Hence each such extrinsic symmetric space is again an orbit $\Ad(K)\cdot \xi$ of the isotropy representation of a `bigger' symmetric space $G/K$, where $\ad(\xi)^3=-\ad(\xi)$.   Obviously, at this point the two following problems occur. Firstly, contrary to the situation of a Euclidean symmetric space we cannot assume without loss of generality that our extrinsic symmetric space $M$ is full and normal. In general, now $M$ cannot be decomposed into a flat and a normal extrinsic symmetric space. Moreover, the smallest affine subspace that contains  $M$ can be degenerate. Hence, our assumptions on $M$ to be full and normal are not adequate. Secondly, Ferus' construction produces pseudo-Riemannian symmetric spaces and only very few of these are semisimple. A list of the occuring semisimple symmetric spaces together with the arising extrinsic symmetric spaces can be found in \cite{Nai}. However, one cannot expect to get such a list in the non-semisimple case. So the main tasks are the following
\begin{enumerate}
\item Find a one-to-one correspondence between arbitrary (i.e., not necessarily full and normal) extrinsic symmetric spaces and suitable infinitesimal objects.
\item Find a practicable method to describe these infinitesimal objects.
\end{enumerate}
The aim of this paper is to give a solution of the first problem. A first step towards a solution can be found in \cite{Kim}, see also~\cite{EK}. Kim generalised Ferus' result and proved that to each full, not necessarily normal extrinsic symmetric space $M\hookrightarrow V$ one can associate the Lie algebra $\fg$ of the transvection group of a bigger symmetric space. However, if $M\hookrightarrow V$ is not normal then it is not an orbit of the isotropy representation of this bigger symmetric space. In this case Kim's result does not allow to reconstruct the extrinsic symmetric space from the associated infinitesimal object. 

Here we will improve Kim's construction and we will show that any full extrinsic symmetric space in a pseudo-Euclidean space is uniquely determined by an object which we will call extrinsic symmetric triple. This object consists of a Lie algebra $\fg$, a $\fg$-invariant non-degenerate inner product on $\fg$ and a pair $(D,\theta)$ consisting of an isometric involution $\theta$ on $\fg$ and a derivation $D$ on $\fg$ such that $D\theta=-\theta D$ and $D^3=-D$. 

To get rid of the inadequate assumption that the extrinsic symmetric space $M$ is full in the ambient Euclidean space we proceed as follows. If $M$ is not full, then we restrict $M$ to the smallest affine subspace $V_M$ that contains $M$. Then $M$ is full in $V_M$. However, the inner product on $V_M$ can be degenerate. Therefore we will consider as an intermediate step also (non-degenerate) full extrinsic symmetric spaces in degenerate inner product spaces. Again we can describe these extrinsic symmetric spaces by an  infinitesimal object, which is now called weak extrinsic symmetric triple. 

We prove that any weak extrinsic symmetric triple is a central extension of an extrinsic symmetric triple. These extensions can be described by a Lie algebra cohomology group. In particular, this solves the extension problem posed in \cite{EK}, see Remark~\ref{REK}.

The above described results of this paper reduce the classification of extrinsic symmetric spaces to the classification of extrinsic symmetric triples and their central extensions. 

We conclude the paper with examples which illustrate the construction of extrinsic symmetric spaces from weak extrinsic symmetric triples. In particular, we will see that non-normal full extrinsic symmetric spaces in pseudo-Euclidean vector spaces do exist. This answers the question from \cite{EK} whether there are full extrinsic symmetric spaces in a pseudo-Euclidean space that are not of Ferus type (in the sense of \cite{EK}).

In the forthcoming paper \cite{Kext2} we will deal with the second problem posed above. We will apply the method of quadratic extensions developed in \cite{KO1, KO2, KOesi} to find a description of the occuring infinitesimal objects by means of a certain cohomology set. Moreover, we will classify Lorentzian extrinsic symmetric spaces explicitly.

I would like to thank Dmitri Alekseevsky,
Jost Eschenburg and Martin Olbrich for many helpful discussions.

\section{Extrinsic symmetric spaces} \label{S2}

In this section we recall some basic facts about extrinsic symmetric spaces. Although we are interested in results about extrinsic symmetric spaces in pseudo-Euclidean affine spaces it will be helpful to consider also extrinsic symmetric spaces in degenerate inner product spaces as an intermediate step.

Let $V$ be a finite-dimensional real affine space. The associated vector space will also be denoted by $V$. Let  $\ip$ be a symmetric bilinear form on the vector space $V$. We will call $\ip$ inner product although we do not  assume that $\ip$ is non-degenerate. The pair $(V,\ip)$ is called inner product space. It is called pseudo-Euclidean if $\ip$ is non-degenerate and Euclidean if $\ip$ is positive definite. An isometry between two inner product spaces $(V_1,\ip_1)$ and $(V_2,\ip_2)$ is an invertible affine map $\ph:V_1\rightarrow V_2$ that satisfies $\ph^*\ip_2=\ip_1$.  It is called linear if it is a linear map from $V_1$ to $V_2$. Then 
$$\grO(V,\ip):=\{A\in\GL(V)\mid \langle Ax,Ay\rangle=\langle x,y\rangle\}\subset \GL(V)$$
is the group of linear isometries from $(V,\ip)$ to itself.
The Lie algebra of this group equals
$$\fo(V,\ip):=\{A\in\fgl(V)\mid \langle Ax,y\rangle+\langle x,Ay\rangle=0\}\subset \fgl(V).$$
The isometry group of $(V,\ip)$ equals $\grO(V,\ip)\ltimes V$, where we understand $V$ as the group of translations of $V$. Often we write $V$ instead of $(V,\ip)$.

In particular, we consider
$V:=\RR^{p,q,r}=\RR^{p+q+r}$ together with the inner product $\ip_{p,q,r}$ given by 
$$\langle x,y\rangle_{p,q,r}=-x_1y_1-\dots-x_py_p+x_{p+1}y_{p+1}+\dots+x_{p+q}y_{p+q}.$$
In the following let $V,V',V_1,V_2,$ etc.~always be inner product spaces.

We will consider submanifolds of $V$. All these submanifolds are assumed to be connected.

Let $M_1\hookrightarrow V_1$ and $M_2\hookrightarrow V_2$ be submanifolds. An extrinsic isometry from $M_1$ to $M_2$ is an 
(affine) isometry $V_1\rightarrow V_2$ that takes $M_1$ to $M_2$. If such an extrinsic isometry from $M_1$ to $M_2$ exists, we will say that $M_1$ and $M_2$ are extrinsic isometric or just that $M_1\hookrightarrow V_1$ and $M_2\hookrightarrow V_2$ are isometric.

In the following we will be interested in a special kind of extrinsic isometries, namely in reflections.  
Let $M\hookrightarrow V$ be a non-degenerate submanifold. For $x\in M$ we denote by $\cT_xM$ the tangent space at $x$ and by $\cN_x M$ the normal space at $x$, where we understand both of these spaces as affine subspaces of $V$. Moreover, let $T_xM$ be the linear subspace of $V$ that is associated with  $\cT_xM$ and put $T^\perp_xM:=(T_xM)^\perp\subset V$. The reflection $s_{ x}$ of $V$ at $\cN _xM$ is the affine map $s_x:V\rightarrow V$ defined by 
\begin{equation}\label{Es}
ds_{x}|_{T_xM}=-\Id,\quad s_{x}|_{{\cal N}_xM}=\Id.
\end{equation} 
\begin{de}\label{Dex}
A non-degenerate submanifold  $M\hookrightarrow V$ is called extrinsic symmetric if $s_x(M)=M$ for all $x\in M$.
\end{de} 
Now we can formulate the classification problem for extrinsic symmetric spaces.\\[2ex]
{\bf Problem\ \ } Determine all extrinsic symmetric spaces in arbitrary pseudo-Euclidean spaces up to extrinsic isometry.\\[2ex]
As said in the introduction the objective of this paper is to reduce this classification problem to a purely algebraic one. Let us start with recalling some basic properties of extrinsic symmetric spaces.

Obviously, every extrinsic symmetric space is an (ordinary) pseudo-Riemannian symmetric space. 

 Let $\nabla$ be the flat connection on $V$. Then $\nabla$ defines 
a connection $\nabla^M$ on $T M$ and a connection $\nabla^\perp$ on $T^\perp M$ 
by 
$$\nabla^M_X Y=\proj_{TM}\nabla_X Y,\quad \nabla^\perp_X 
N=\proj_{TM^\perp}\nabla_X N$$ 
for vector fields $X$ and $Y$ on $M$ and sections $N$ in $T^\perp M$, 
where $\proj_{TM}$ and $\proj_{TM^\perp}$ denote the orthogonal projections of $V$ onto $TM$ and $TM^\perp$, respectively. 

It is easy to see that the second fundamental form $\alpha$ of an extrinsic symmetric space $M\hookrightarrow V$ is parallel, i.e., $\nabla^\perp\alpha=0$. Conversely, a complete non-degenerate submanifold $M\hookrightarrow V$ with parallel second fundamental form is extrinsic symmetric \cite{S}.

Now let $M\hookrightarrow V$ be extrinsic symmetric. 
Let $\gamma$ be a geodesic in $M$. An extrinsic isometry $\Phi$ 
of $M$ is called transvection along~$\gamma$ by $c\in\RR$ if $\Phi$ leaves invariant 
$\gamma$ and if $d\Phi|_{\gamma(t)}:T_{\gamma(t)}M\rightarrow 
T_{\gamma(t+c)}M$ 
and $d\Phi|_{\gamma(t)}:T^\perp_{\gamma(t)}M\rightarrow T^\perp_{\gamma(t+c)}M$ 
are parallel displacements along~$\gamma$. Similar to the case of ordinary symmetric spaces one proves that the subgroup   
$$K:=\langle s_x\circ s_y\mid x,y\in M\rangle.$$ 
of the group of extrinsic isometries of $M\hookrightarrow V$ contains all transvections of $M$ and is generated by them.
Therefore $K$ is called transvection group of the extrinsic symmetric space $M\hookrightarrow V$.  If $V$ is definite, then the above defined transvection 
group is isomorphic to the (ordinary) transvection group of the (ordinary) pseudo-Riemannian 
symmetric space $M$. For pseudo-Riemannian or degenerate $V$ this is not true in general as Example~\ref{Ex*} will show.

The transvection group $K$ acts transitively on $M$. Let us fix a point $x_0\in 
M$ and let $K_+\subset K$ be the stabilizer of $x_0$. Then $M=K/K_+$. The 
reflection $s_{x_0}$ defines an involution $s$ on $K$ by $s(k)=s_{x_0}ks_{x_0}$ 
for $k\in K$. Let $\fk$ be the Lie algebra of $K$. The involution $s$ on $K$ induces an involution on $\fk$. Let us introduce the following convention, which we will use throughout this paper.\\[1ex]
{\bf Notation.}
Given a Lie algebra $(\fg,\theta)$ with involution we will denote by 
$\fg_+$ and $\fg_-$ the eigenspaces of $\theta$ with eigenvalue $1$ and $-1$, 
respectively.\\[1ex]
For the Lie algebra $\fk$ of the transvection group and its involution $s$ defined above one can prove   $[\fk_-,\fk_-]=\fk_+$. The subalgebra $\fk_+$ is 
the Lie algebra of $K_+$. Moreover, $\fk_-$ can be identified with the vector 
space $V_-:=T_{x_0}M$.  More exactly, given an element $X\in\fk_-$ we consider 
the one-parameter group $\Phi_t:=\exp tX$ of extrinsic isometries of $M$. These extrinsic isometries are exactly the transvections along the geodesic $\gamma(t)=(\exp tX)(x_0)$. We 
identify 
$X$ with the tangent vector of $\gamma $ in $x_0$. Let 
\begin{equation}\label{dT}
T: V_- \longrightarrow \fk_-,\quad u\longmapsto T_u
\end{equation}
be the inverse of this map. 
We understand each $T_u$ as a map from the affine space to the vector space $V$.
In particular, 
\begin{equation}\label{Tuu}
T_u(x_0)=u.
\end{equation}
The map $T$ is related to the second fundamental 
form $\alpha$ of $M$ and to the shape operator $A$ of $M$ at the point $x$. 
Namely, the differential $dT_u\in\fo(V,\ip)$ of the map $T_u$ satisfies 
\begin{equation}\label{ET} 
dT_u(v)=\alpha(u,v),\quad dT_u(\eta)=-A_\eta u 
\end{equation} 
for all $u,v\in V_-=T_{x_0}M$ and $\eta \in V_+:=T^\perp_{x_0}M$. The proof in the case of a Euclidean space $V$ given in \cite{EH} remains valid for arbitrary inner product spaces as long as $M\hookrightarrow V$ is non-degenerate. Since these equations are essential for our considerations in the following sections we will recall the argument here.

For all $x\in V$ we identify the tangent space $T_xV$ with $V$. Take $X\in T_{x_0} V\cong V$. Then we have 
\begin{eqnarray*}
dT_u(X)&=& \frac d{dt} T_u(x_0+tX)|_{t=0}=\ \frac d{dt}\Big(\frac d{ds}\Phi_s\Big)|_{s=0}(x_0+tX)|_{t=0}\\
&=&\frac d{ds}\frac d{dt}\Phi_s(x_0+tX)|_{t=0}|_{s=0}\ =\ \frac d{ds} d\Phi_s(X)|_{s=0}
\ =\ \nabla_u(d\Phi_s(X)).
\end{eqnarray*}
Moreover, since $d\Phi_s$ restricted to $TM|_\gamma$ and $T^\perp M|_\gamma$ is the parallel translation we have
\begin{eqnarray*}
\nabla_u(d\Phi_s(v))&=&\nabla_u^M(d\Phi_s(v))+\alpha(u,d\Phi_0(v))\ =\ \alpha(u,v)\\
\nabla_u(d\Phi_s(\eta))&=&-A_{d\Phi_0(\eta)}u+\nabla^\perp_u(d\Phi_s(\eta))\ =\ -A_\eta u
\end{eqnarray*}
for all $v\in V_-$ and $\eta\in V_+$, which proves~(\ref{ET}).

Let $R^M$ and $R^\perp$ be the curvature of $\nabla^M$ and $\nabla^\perp$, respectively. Similar to the case of ordinary symmetric spaces one proves
\begin{equation}\label{R}
R^M(u,v)w=-[T_u,T_v](w),\quad R^\perp(u,v)\eta=-[T_u,T_v](\eta)
\end{equation}
for $u,v,w\in V_-$ and $\eta\in V_+$.

\section{Full submanifolds} \label{S3}
Let $V$ be an inner product space.
Recall that a submanifold $M\hookrightarrow V$ is called full if $M$ is not contained in a proper affine subspace of $V$. The following consideration will show that in order to classify extrinsic symmetric spaces in pseudo-Euclidean vector spaces (i.e., in non-degenerate inner product spaces) it is sufficient to classify full extrinsic symmetric spaces in (arbitrary) inner product spaces.

If $M\hookrightarrow V'$ is extrinsic symmetric  and if $V\subset V'$ is an affine subspace that contains $M$, then  $M\hookrightarrow V$ is again extrinsic symmetric. Moreover, restricting our considerations to $V$ we don't loose any information as the following proposition will show.  
\begin{pr}\label{2.1.}
 Let $f:M\hookrightarrow V$ be an extrinsic symmetric space and let $V'$ be a further inner product space.
\begin{enumerate} 
\item If $i:V\rightarrow V'$ is an affine isometric embedding, then $i\circ f:M\hookrightarrow V'$ is also extrinsic symmetric.
\item If $i_1:V\rightarrow V'$ and $i_2:V\rightarrow V'$ are affine isometric embeddings into a non-degenerate space $V'$, then $i_1\circ f:M\hookrightarrow V'$ and $i_2\circ f:M\hookrightarrow V'$ are isometric.
\end{enumerate}
\end{pr}
\proof 
ad 1.) We identify $M$ with its image in $V$ and $V$ with its image in $V'$. Take $x\in M$ and let $s_x'$ be the corresponding reflection in $V'$ and $s_x$ the reflection in $V$. Since
$$V=x+T_xM\oplus(T_xM)^{\perp_V}=x+T_xM\oplus \Big((T_xM)^{\perp_{V'}}\cap V\Big)$$
we have $s'_x(V)\subset V$ and $s'_x|_V=s_x:V\rightarrow V$. In particular, $s_x'(M)=M$.

ad 2.) Since translations are always isometries we may restrict ourselves to the case, where $i_1$ and $i_2$ are linear.
The map
$$i_1(V)\stackrel{(i_1)^{-1}}{\longrightarrow} V\stackrel{i_2}{\longrightarrow}V'$$
is an isometric embedding. Since $V'$ is non-degenerate it can be extended to an isometry of $V'$ by Witt's theorem. This isometry maps $i_1\circ f(M)$ onto $i_2\circ f(M)$.
\qed

Note that the second assertion is not true for degenerate $V'$. Furthermore, if $M\hookrightarrow V$ is a submanifold and $V_1$ and $V_2$ are affine subspaces of $V$ that contain $M$, then $M\hookrightarrow V_1$ and $M\hookrightarrow V_2$ are not necessarily isometric even if $V_1$ and $V_2$ have the same signature. 
To avoid difficulties we will often consider $M$ as a submanifold of a special subspace of $V$, which is uniquely determined by $M$, namely as a submanifold of the smallest subspace of $V$ that contains $M$. Let us study this subspace in more detail.

For a  non-degenerate submanifold $M\hookrightarrow V$ we define the vector space
$$W_M:=  \sum_{x\in M}T_xM\subset V.$$
The following statement can be proved easily.
\begin{lm} \label{2.2.}
Let $M\hookrightarrow V$ be a non-degenerate submanifold and let $x_0\in M$ be fixed.
The smallest affine  subspace that contains $M$ is equal to 
$$V_M:=x_0+W_M.$$ In particular, $M\hookrightarrow V$ is full if and only if $V=W_M$. 
\end{lm}

\begin{pr}\label{Pfull} If $M\hookrightarrow V$ is an extrinsic symmetric space and $x_0\in M$, then $$W_M= \im \alpha_{x_0}\oplus V_-.$$
\end{pr}
\proof Let us first prove the following local version. Let $U$ be an open set in $M$ that is a normal neighbourhood of each of its points and take $x\in U$. In particular,  there is a geodesic in $U$ from $x$ to $y$ for all $y\in U$. Then  $$W_{U}:=\sum_{y\in U} T_yM =\im \alpha_{x} +T_xM.$$
The inclusion $T_xM\subset W_{U}$ is obvious. Let us verify that also $\im \alpha_{x}\subset W_{U}$ holds. Take $w\in \im \alpha_{x}$. Then $w= dT_u(v)$ for suitable $u,v\in T_xM$, where $T$ is defined as in (\ref{dT}) for $x$ instead of $x_0$. Then $\exp t T_u$ is a one-parameter group of extrinsic isometries of $M$. Hence  $\mu(t):d(\exp (t T_u))(v) =\exp( t dT_u) (v)\in T_{\exp tT_u(x)}M$ is a curve is $W_{U}$. Since $W_{U}$ is a vector space also $\mu'(0)=dT_u (v)=w$ is in $W_{U}$. 

Now let us show that $W_{U}\subset \im \alpha_{x} +T_xM$ holds. Take $w\in T_yM$ for $y\in U$. Let $\gamma$ be a geodesic with $\gamma(0)=x$ and $\gamma(1)=y$. Let $\Phi$ be the transvection along $\gamma$ that sends $x$ to $y$. Then we have $d\Phi = d(\exp T_u)=\exp dT_u$ for $u=\gamma'(0)$. In particular, $w=\exp dT_u (v)$ for some $v\in T_xM$. Hence 
$$w=v+dT_u(v)+\frac 12 (dT_u)^2(v)+\frac 1{3!}(dT_u)^3(v)+\dots.$$
Since by (\ref{ET}) the image of $dT_u$ is contained in $\im \alpha_x\oplus T_xM$ this implies $w\in  \im \alpha_{x} +T_xM$.

Now let us prove the assertion of the proposition. Let ${\cal U}=\{U_i\mid i\in I\}$ be  a covering of $M$ by open sets that are normal with respect to each of their points. Fix a $U_0\in{\cal U}$ with $x_0\in U_0$ and define $I_1=\{ i\in I\mid W_{U_i}=W_{U_0}\}$, $I_2=\{ i\in I\mid W_{U_i}\not=W_{U_0}\}$. Then $\bigcup_{i\in I_1}U_i$ and  $\bigcup_{i\in I_2}U_i$ are disjoint. Indeed, by the above considerations  $x\in U_{i}\cap U_{j}$ for $i\in I_1$, $j\in I_2$, would imply $W_{U_0}=W_{U_i}=W_{U_j}\not=W_{U_0}$, a contradiction. Since $M$ is connected we obtain $M=\bigcup_{i\in I_1}U_i$. Hence $W_M=\sum_{i\in I_1} W_{U_i}= W_{U_0}=\im \alpha_{x_0}\oplus V_-.$\qed

\section{Extrinsic symmetric triples and associated extrinsic \\ symmetric spaces} 
\label{S4}

It is well known  that there is a one-to-one correspondence between simply connected pseudo-Riemannian symmetric spaces and certain algebraic objects which we call symmetric triples. The aim of the next two sections is to establish a similar correspondence for full extrinsic symmetric spaces.

Let us first recall the notion of a metric Lie algebra. A metric Lie algebra is a pair $(\fg,\ip)$ consisting of a Lie algebra $\fg$ and a non-degenerate inner product $\ip$ that is invariant under the adjoint representation of $\fg$. A metric Lie algebra with involution  $(\fg,\ip,\theta)$ consists of a metric Lie algebra $(\fg,\ip)$ and an 
isometric involution $\theta$. Such a  metric Lie algebra with involution is called symmetric triple if the eigenspace decomposition $\fg=\fg_+\oplus\fg_-$ with respect to $\theta$ satisfies $\fg_+=[\fg_-,\fg_-]$. 

Let $(\fg,\ip,\theta)$ be a symmetric triple. Let $D$ be a derivation of 
$\fg$ that satisfies 
$D^3=-D$ and $D\circ \theta=-\theta\circ D$. In particular, $D$ is semisimple and the eigenvalues of $D$ are in 
$\{i,-i,0\}$. Let us define 
$$\fg^+:=\Ker D,\quad \fg^-:=\{X\in\fg\mid D^2(X)=-X\}.$$ Then 
$$\tau_D:\fg\longrightarrow\fg,\quad \tau_D|_{\fg^+}=\Id,\ 
\tau_D|_{\fg^-}=-\Id$$ 
defines an involution $\tau_D$ on $\fg$. The subspaces $\fg_+$ and $\fg_-$ are invariant under $\tau_D$. In particular, $\tau_D$ and 
$\theta$ commute.
We introduce the notation 
\begin{equation}\label{not} 
\fg_+^+:=\fg_+\cap\fg^+,\ \fg_+^-:=\fg_+\cap\fg^-,\ 
\fg_-^+:=\fg_-\cap\fg^+,\ \fg_-^-:=\fg_-\cap\fg^-. 
\end{equation} 

\begin{de}\label{Dq} 
An extrinsic symmetric triple $(\fg,\ip,\Phi)$ consists of a metric Lie 
algebra $(\fg,\ip)$ and a pair $\Phi=(D,\theta)$, where 
$\theta\in\Aut(\fg)$ is an isometric involution, and $D\in\fso(\fg)$ is an anti-symmetric derivation 
satisfying $D^3=-D$  and $
[\fg_+^-,\fg_+^-]=\fg^+_+.$
\end{de} 
In this situation, $(\fg_+,\ip|_{\fg_+\times\fg_+},\tau_D|_{\fg_+})$ is a symmetric triple.
\begin{de}
We want to weaken the above definition by supposing non-degeneracy of $\ip$ only on $\fg_+\oplus\fg_-^-$. We then speak of a weak extrinsic symmetric triple.
\end{de}
It is easy to prove the following statement. 
\begin{lm} \label{Leta} 
If $(\fg,\ip,\Phi)$  is a (weak) extrinsic symmetric triple with $\Phi=(D,\theta)$, then $\, 
D|_{\fg_+^-}: 
\fg_+^-\rightarrow \fg_-^-$ is an isometry. 
\end{lm} 
\begin{de} 
Two (weak) extrinsic symmetric triples $(\fg_i,\ip_i,\Phi_{i})$, $i=1,2$, with 
$\Phi_{i}=(D_i,\theta_i)$ 
are called isomorphic if there is an isomorphism $h:\fg_1\rightarrow \fg_2$ 
of Lie algebras such that 
$$\ip_1=h^*\ip_2,\quad D_2=h 
D_1h^{-1},\quad \theta_2=h\theta_1h^{-1}.$$ 
\end{de} 
Let $(\fg,\ip,\Phi)$ be a weak extrinsic symmetric triple with $\Phi=(D,\theta)$. 
We consider $\fg_-$ together with the restriction of $\ip$ to $\fg_-$ as an inner product space. The isometry group of $\fg_-$ equals $\Iso(\fg_-)=\grO(\fg_-)\ltimes \fg_-$.  Let $\phi$ be the map from $\fg_+$ to the Lie algebra $\frak{iso}(\fg_-)=\fso(\fg_-)\ltimes\fg_-$ of $\Iso(\fg_-)$ defined by
$$\phi(X)=((\ad X)|_{\fg_-},-D(X)),\quad X\in\fg_+.$$
Since $D$ is a derivation $\phi$ is a homomorphism. Moreover, $\phi$ is injective. Indeed, assume $\phi(X) =0$. Then $\ad(X)|_{\fg_-}=0$ and $DX=0$. Hence $[X,Du]=D[X,u]=0$ for all $u\in\fg_-^-$. Thus $[X,\fg_+^-]=0$, which implies $X=0$ since $X\in\fg_+^+$ and $(\fg_+,\ip|_{\fg_+\times\fg_+},\tau_D|_{\fg_+})$ is a symmetric triple. Hence $\ker\phi=0$. In particular, we may identify the Lie algebra of $G_+$ with $\fg_+$.
Let $G_+\subset \Iso(\fg_-)$ be the group 
$$G_+=\langle\, \exp(\phi(X))\mid X\in\fg_+\,\rangle.$$ 
We define 
$$M_{\fg,\Phi}:=G_+(0)\subset\fg_-.$$ 

\begin{pr} \label{POrbit} For any weak extrinsic symmetric triple $(\fg,\ip,\Phi)$ the manifold
$M:=M_{\fg,\Phi}$ is an extrinsic symmetric 
space in $\fg_-$ whose transvection group equals~$G_+$. The isometry $T$ defined in (\ref{dT}) equals
$T=D|_{\fg_-^-}: \fg_-^-\rightarrow \fg_+^-$.

The (ordinary) symmetric 
space $M$ with the induced metric is associated with the symmetric triple 
$(\fg_+,\ip|_{\fg_+\times\fg_+},\tau_D|_{\fg_+})$. 
\end{pr} 
\proof 
The group $G_+$ acts by extrinsic isometries on $M$.  Let $G_+^+\subset G_+$ be 
the stabilizer of $0\in M$. Now we consider the involutions $\theta$ and 
$\tau_D$ and the corresponding eigenspace decomposition of $\fg$. 
Then $\fg_+^+$ is the Lie algebra of $G_+^+$ and we have 
\begin{eqnarray}\label{TT} 
T_0 M&=&\Big\{\textstyle{\frac d{dt}} \exp(\phi(tX))(0)|_{t=0} \ \Big| 
X\in \fg_+\Big\} \nonumber\\
&=&\{D(X)\mid X\in \fg_+\}=\fg_-^-,\nonumber \\
 T^\perp_0 M&=&\fg_{-}\cap(\fg_-^-)^\perp=\fg_-^+. 
\end{eqnarray} 
Now we prove that $M$ is extrinsic symmetric in $\fg_-$. First we consider 
$0\in M$ and define $s_0=\tau_D|_{\fg_-}$. Because of (\ref{TT}) we have $ds_0=-\Id$ on $T_0 M$ 
and $s_0=\Id$ on $\cN_0M=T^\perp_0 M$. Hence $s_0$ is the reflection at $\cN_0M$. Since $\tau_D$ is an 
isomorphism of $\fg$ commuting with $D$ we have 
$$\tau_D\circ  \exp((\ad X)|_{\fg_-},-D(X))= \exp((\ad \tau_D(X))|_{\fg_-},-D(\tau_D (X)))\circ \tau_D,$$
hence
$$\tau_D\circ  \exp(\phi(X))= \exp(\phi( \tau_D(X)))\circ \tau_D.$$
Since $\tau_D(X)\in\fg_+$ for all $X\in\fg_+$ we obtain 
$\tau_D G_+=G_+\tau_D$. 
We get 
$$\tau_D(G_+(0))= 
G_+\cdot \tau_D (0)= G_+(0)=M,$$ 
hence $s_0(M)=M$.  Hence $s_0$ satisfies the conditions 
in Definition \ref{Dex}. Now let $x\in M$ be arbitrary. Then $x=g(0)$ for 
some $g\in G_+$ and $s_x:=g s_0 g^{-1}$ is the reflection at $\cN_xM$. Hence, 
$M$ is an extrinsic symmetric space. 

The assertion on $T$ follows from 
$$T^{-1}(A)=A(0)=-D(A)$$
for $A\in \fk_-=\fg_+^-$ and from $(D|_{\fg^-})^2=-\Id$.

We have seen that $M\in M_{\fg,\Phi}$ is the homogeneous space $G_+/G_+^+$ and that $\fg_+$ is the Lie algebra of $G_+$ and $\fg_+^+$ is the Lie algebra of $G_+^+$. In order to prove that $M$ is associated with the symmetric triple 
$(\fg_+,\ip|_{\fg_+\times\fg_+},\tau_D|_{\fg_+})$ it remains to prove that 
the metric on $M$ which is induced by the inner product 
$\ip|_{\fg_-\times\fg_-}$ corresponds to $\ip|_{\fg_+\times\fg_+}$. But this 
follows from the fact that $D:\fg_+^-\rightarrow \fg_-^-= T_0 M$ is an isometry. 

It remains to show that $G_+$ is the transvection  group of $M\hookrightarrow V$. Since $G_+$ is a connected subgroup of $\Iso(\fg_-)$ it suffices to show that its 
Lie algebra coincides with the Lie algebra $\fk$
 of the transvection group of $M\hookrightarrow V$. Using $G_+\tau_D=\tau_D G_+$ we see that 
$G_+$ contains $s_xs_y$ for all $x,y\in M$ and therefore also the whole transvection group. Thus $\fk\subset\fg_+$. On the other hand, there is a surjective homomorphism from $\fk$ onto $\fg_+$ since by the above considerations $\fg_+$ is the Lie algebra of the transvection group of the symmetric space $M$. 
Hence $\fk=\fg_+$.
\qed 

\begin{re}{\rm
Note that we here are in a special situation, where the Lie algebras of the ordinary transvection group and the extrinsic transvection group coincide. In general, the transvection group of an extrinsic symmetric space $M\hookrightarrow V$ can be much larger than the  transvection group of the symmetric space $M$ if $V$ is pseudo-Euclidean.  This phenomenon will be illustrated by Example~\ref{Ex*}.  
}
\end{re}
\begin{re}\label{R4.2.}
 {\rm If $D=\ad (\xi)$ is an inner derivation (and if $\fg$ is non-degenerate), then we are in the situation described in the introduction. There we associated with $\fg=\fg_+\oplus\fg_-$ and $\xi$ the extrinsic symmetric space $M_{\fg,\xi}:= \langle \exp ((\ad X)|_{\fg_-})\,|\, X\in\fg_+\rangle\cdot \xi\subset\nolinebreak\fg_-$. Let us compare this space with $M_{\fg,\Phi}$. Since $D=\ad (\xi)$ we get $\exp(\phi(X))(0)=\exp((\ad X)|_{\fg_-})(\xi)-\xi$, hence 
$M_{\fg,\Phi}$ differs from $M_{\fg,\xi}$ by a translation by $-\xi$.
}
\end{re}

\begin{pr} \label{Propiso}
Let $(\fg_i,\ip_i,\Phi_i)$, $i=1,2$ be weak extrinsic symmetric triples and let 
$M_i= M_{\fg_i,\Phi_i}$ be the associated 
extrinsic symmetric spaces in $(\fg_i)_-$. Then $M_1\hookrightarrow(\fg_1)_-$ and $M_2\hookrightarrow(\fg_2)_-$ are isometric 
if and only if $(\fg_1,\ip_1,\Phi_1)$ and 
$(\fg_2,\ip_2,\Phi_2)$ are isomorphic. 
\end{pr} 
\proof 
Assume that there is an isometry 
$f:(\fg_1)_-\rightarrow (\fg_2)_-$ such that 
$f(M_1)=M_2$.  Since $0\in M_1$ and $0\in M_2$ and $M_2$ is 
extrinsic symmetric we may assume that 
$f(0)=0$.  Now we consider $f$ as an isometry $f:M_1\rightarrow M_2$. Then conjugation by $f$ defines an isomorphism $F$ from the (ordinary) transvection group of $M_1$ to the transvection group of $M_2$ and $dF$ is an isomorphism of the associated symmetric triples. According to Prop.~\ref{POrbit} $dF$ can be considered as an isomorphism from  
$((\fg_1)_+,\ip_1|,\tau_{D_1}|_{(\fg_1)_+})$ to $((\fg_2)_+,\ip_2|,\tau_{D_2}|_{(\fg_2)_+})$. Moreover, $dF(T_u)=T_{f(u)}$. 

Now we define $h:=dF\oplus f:(\fg_1)_+\oplus (\fg_1)_-\rightarrow 
(\fg_2)_+\oplus (\fg_2)_-$. 
By definition we have 
$h((\fg_1)_+)=(\fg_2)_+$ and $h((\fg_1)_-)=(\fg_2)_-$, hence 
$h\theta_1h^{-1}=\theta_2$. Moreover, $h:\fg_1\rightarrow \fg_2$ is an isometry since $dF$ and $f$ are so. We show that $h$ is also an isomorphism of Lie algebras. For this it remains to prove  that $h[A,X]=[hA,hX]$ for all $A\in (\fg_1)_+$ and all $X\in  (\fg_1)_-$ and that $h[X,Y]=[hX,hY]$ for all $X,Y\in(\fg_1)_-$. For $A=T_u\in  (\fg_1)_+^-$ the first equation follows from (\ref{ET}) and the fact that $f$ is an extrinsic isometry. From $[(\fg_1)_+^-,(\fg_1)_+^-]=(\fg_1)_+^+$ now it follows that the first equation is also true for  $A\in  (\fg_1)_+^+$. Let us now verify the second equation. Take $X,Y\in(\fg_1)_-$. Then 
\begin{eqnarray*}
\langle  h[X,Y], hA\rangle_2 &=& \langle [X,Y],A\rangle_1 \,=\, -\langle X,[A,Y]\rangle_1\,= \,-\langle hX,h([A,Y])\rangle_2\\ &=&-\langle hX, [hA,hY]\rangle_2 \ =\ \langle [hX,hY], hA\rangle_2
\end{eqnarray*}
for all $A\in(\fg_2)_+$, hence $h[X,Y]=[hX,hY]$.

It remains to prove that $hD_1(Z)=D_2h(Z)$ holds for all $Z\in\fg_1$. Again we start with the case $Z=T_u\in(\fg_1)_+^-$ and see  
$$hD_1(T_u)=h(D_1)^2(u)=-h(u)=-f(u)=(D_2)^2(f(u))=D_2T_{f(u)}=D_2h(T_u).$$ Because of $[(\fg_1)_+^-,(\fg_1)_+^-]=(\fg_1)_+^+$ this implies that the assertion is also true for all $Z\in(\fg_1)_+^+$. Now take $Z\in(\fg_1)_-$. Then 
\begin{eqnarray*}
\langle hD_1(Z),A\rangle_2 &=& \langle D_1(Z), h^{-1}A\rangle_1 \ =\  -\langle Z,D_1h^{-1}A\rangle_1\ =\ -\langle h(Z),hD_1h^{-1}A\rangle_2\ \\ &=& -\langle h(Z),D_2(A)\rangle_2 \ = \  
\langle D_2 h(Z),A\rangle_2
\end{eqnarray*}
for all $A\in (\fg_2)_+$, hence $hD_1(Z)=D_2h(Z)$.

Vice versa, let $h:\fg_1\rightarrow\fg_2$ be an isomorphism such that $h^*\ip_2=\ip_1$, $h\circ\theta_1=\theta_2\circ h$ and $h\circ D_1=D_2\circ h$.  Then $h$ restricts to an extrinsic isometry from $M_1\hookrightarrow (\fg_1)_-$ to $M_2\hookrightarrow (\fg_2)_-$. \qed

Since we are especially interested in full extrinsic symmetric spaces we now want to clarify which weak extrinsic symmetric triples are associated with these spaces.
\begin{de} A weak extrinsic symmetric triple $(\fg,\ip,\Phi)$ is called full if it satisfies $[\fg_+^-,\fg_-^-] =\fg_-^+$ or, equivalently, $[\fg^-,\fg^-]=\fg^+$. 

\end{de}
\begin{co} Let $(\fg,\ip,\Phi)$ be a weak extrinsic symmetric triple and let $M\in M_{\fg,\Phi}$ be an associated extrinsic symmetric space in $\fg_-$. Then $M\hookrightarrow \fg_-$ is full if and only if $(\fg,\ip,\Phi)$ is full. 
\end{co}
\proof {}From the definition of the $G_+$-action on $\fg_-$ we know that
$$dT_u(v)=\ad(T_u)|_{\fg_-}(v)=[T_u,v]$$
holds for all $u,v\in\fg_-^-$. By (\ref{ET}) this implies $\Im\alpha_0=[\fg_+^-,\fg_-^-]$. Now the assertion follows from Prop.~\ref{Pfull}. \qed
\section{The inverse construction} 
In this section we will study a construction that may be 
considered as a converse of Proposition \ref{POrbit}. Given a full extrinsic 
symmetric space $M\subset V$ of an inner product space 
this construction yields a weak extrinsic symmetric 
triple $(\fg,\ip,\Phi)$ such that $M=M_{\fg,\Phi}$. This construction is based on an idea due to Ferus. In \cite{F2,F3} Ferus studied extrinsic symmetric spaces in Euclidean ambient spaces. To each full and compact extrinsic symmetric space $M\hookrightarrow \RR^n$ Ferus associates a compact semisimple symmetric triple $(\fg,\ip,\theta)$ and an element $\xi\in\fg_-$ such that $M$ is extrinsic isometric to the orbit $G_+\cdot\xi\hookrightarrow \fg_-$, where $G_+=\langle \exp((\ad X)|_{\fg_-})\mid X\in\fg_+\rangle$. The element $\xi\in\fg_-$ has the property $\ad(\xi)^3=-\ad(\xi)$. Ferus uses the language of Jordan triple systems. A more elementary 
version of this construction was given by Eschenburg and Heintze in \cite{EH}.  Naitoh  \cite{Nai} generalised this construction to the case of pseudo-Riemannian extrinsic symmetric spaces in (non-degenerate) pseudo-Euclidean spaces satisfying two additional conditions. One of these conditions is equivalent to the fullness condition. The other one is equivalent to the condition to be  normal introduced in \cite{KOesi}. A non-degenerate submanifold $M\hookrightarrow V$ is called normal if $$\bigcap_{x\in M} \cN_xM\not=\emptyset.$$ 
In \cite{EK} Eschenburg and Kim give a sufficient condition for an extrinsic symmetric space to be normal in terms of the shape operator $A_h$ with respect to the mean curvature vector $h$. They prove that a full extrinsic symmetric space in a pseudo-Euclidean space is normal unless $A_h^2=0$. However, there are full and normal extrinsic symmetric spaces for which $A_h^2=0$ holds, see, e.g., Example~\ref{Ex3}. 

In \cite{Kim} J.\,R.\,Kim studies a generalisation of the 
construction to full but non-normal pseudo-Riemannian extrinsic symmetric spaces in pseudo-Euclidean spaces. He associates a metric Lie algebra with involution $(\fg,\ip,\theta)$ to each of these spaces. However, in this case it is not possible to find an element $\xi\in\fg_-$ such that $M=G_+\cdot\xi$. 

Here we will extend Kim's construction. Firstly, we will generalise the construction to full extrinsic symmetric spaces is arbitrary inner product spaces.  Secondly we will find a (not necessarily inner) derivation $D$ such that $M$ can be reconstructed from $(\fg,\ip,\theta)$ and $D$. We will see that $(\fg,\ip,(D,\theta))$ is a weak extrinsic symmetric triple.

Let us now describe the construction. Let $(V,\ip_V)$ be an inner product space  and let $M\hookrightarrow V$ be an extrinsic symmetric space. Suppose that $M$ is full. We fix an origin $0\in V$  such that $0\in M$. Using again the notation $V_-=T_{0}M$, $V_+=T_{0}^\perp M$ introduced in Section~\ref{S2} we can write $V=V_+\oplus V_-$.  Since $M$ is full $\fk$ is isomorphic to the Lie algebra of the transvection group of the symmetric space $M$. Indeed, let $k:V\rightarrow V$ be a transvection of $M\hookrightarrow V$ and assume that $k|_M$ is the identity. Then $dk|_{T_xM}$ is the identity for all $x\in M$, hence $dk|W_M=\Id$. Thus $k$ is the identity on $V=V_M$. The fact that both transvection groups coincide will be important  for the following construction.

Besides the Lie algebra $\fk$ of the transvection group there is a larger Lie algebra associated with $M\hookrightarrow V$, which we will describe now. As a vector space it equals 
$$\fg:=\fk\oplus V.$$  
On $\fk_-$ we choose the inner product in such a way that $T:V_- \rightarrow 
\fk_-$ is an isometry. Since $\fk$ is isomorphic to the Lie algebra of the transvection group of the symmetric space $M$ we can extend this inner product on $\fk_-$ in a unique 
way to a non-degenerate inner product $\ip_\fk$ on $\fk$ that is invariant 
under $\fk$ and $s$ \cite{CP}. Now we define an inner product on 
$\fg=\fk\oplus V$ by $\ip=\ip_\fk\oplus \ip_V$. 

On $\fg$ we consider the antisymmetric bilinear map 
$\lb:\fg\times\fg\rightarrow\fg$ defined by 
\begin{itemize} 
\item[(i)] $\lb$ restricted to ${\fk\times\fk}$ equals the Lie bracket of 
$\fk$, 
\item[(ii)] $[A,v]=dA(v)$ for $A\in\fk$ und $v\in 
V$, 
\item [(iii)] $\lb$ restricted to $V\times V$ is given by the conditions $[V,V]\subset \fk$ and  
$\langle A,[x,y] \rangle=\langle [A,x],y\rangle$ for all 
$A\in\fk$ and $x,y\in V$. This is a correct definition since $\ip_\fk$ is non-degenerate.
\end{itemize} 
Moreover, we have an involution $\theta$ on $\fg$ given by $\fg_+=\fk$, 
$\fg_-=V$ and a map $D$ with $D^3=-D$ given by
$$ D|_{\fk_-}=-T^{-1},\quad D|_{V_-}=T,\quad D|_{\fk_+\oplus V_+}=0. $$

The next theorem will show that $(\fg,\ip,(D,\theta))$ is a weak extrinsic symmetric triple. Note that if $V$ is non-degenerate then $\fg$ is non-degenerate and if $M$ is normal, then $D$ is an inner derivation, cf~Remark~\ref{fe}. In this case our construction yields the same object as the constructions due to Ferus and Naitoh. Hence, for full and  normal extrinsic symmetric spaces in pseudo-Euclidean spaces the theorem is known.  The proof of the Jacobi identity is the key point of the theorem. We will verify Jacobi identity for all combinations of three vectors from $V_+$, $V_-$, $\fk_+$ and $\fk_-$. In the case where not all three vectors are from $V_+$ the proof of Jacobi identity is essentially the same as in the above cited papers. Nevertheless, we will repeat this part of the proof to convince the reader that all remains true even if $V_+$ is degenerate.  The interesting part is the case where all three vectors are from $V_+$.  Here neither the original proof due to Ferus nor their pseudo-Euclidean versions can be translated to arbitrary inner product spaces since they make essential use of the assumption that $V$ is non-degenerate.

\begin{theo} Let $M\subset V
$ be a  full extrinsic 
symmetric space. Let $\fg,\ip$ and $\Phi=(D,\theta)$ be as constructed 
above. 
Then $(\fg,\ip,\Phi)$ is a weak extrinsic symmetric triple and $M=M_{\fg,\Phi}$. It is an extrinsic symmetric triple, if $V$ is non-degenerate.
\end{theo} 
\proof Since $M\hookrightarrow V$ is non-degenerate the inner product $\ip$ is non-degenerate on $\fk\oplus V_-$. The invariance of $\ip$ with respect to $\lb$ can be checked easily. It uses the $\fk$-invariance of $\ip_\fk$, the definition of the bracket and that $dA$ is contained in $\fo(V,\ip_V)$ for each $A\in\fk$.

Equation (\ref{ET}) gives $[\fk_-,V_+]\subset V_-$ and $[\fk_-,V_-]\subset V_+$. Moreover, $[\fk_+,V_+]\subset V_+$ and $[\fk_+,V_-]\subset V_-$ since $\fk_+$ is the Lie algebra of the stabilizer of $0$. By invariance of $\ip$ these equations imply $\langle [V_+,V_+],\fk_-\rangle=\langle [V_-,V_-],\fk_-\rangle=0$ and $\langle [V_+,V_-],\fk_+\rangle=0$. Since $\ip$ is non-degenerate on $\fk$ we obtain 
$$ [V_+,V_+]+[V_-,V_-]\subset \fk_+, \quad [V_+,V_-]\subset \fk_-.$$
Let us prove the Jacobi identity for $\lb$.
We define $\omega:\fg\times\fg\times\fg \longrightarrow \fg$ by
$$\omega(X_1,X_2,X_3)=[[X_1,X_2],X_3]+[[X_2,X_3],X_1]+[[X_3,X_1],X_2]$$
Then $\omega|_{\fk\times\fk\times\fk}=0$ since $\fk$ is a Lie algebra. Furthermore, $\omega|_{\fk\times\fk\times V}=0$. Indeed, 
$$\fk\longrightarrow \fo(V,\ip),\quad A\longmapsto dA$$ is a Lie algebra homomorphism, hence 
$$[[A,B],v]=d[A,B](v)=(dAdB-dBdA)(v)=[A,[B,v]]-[B,[A,v]].$$
Now we consider $\omega_{\fk\times V\times V}$. Since $\omega(\fk,V,V)\subset \fk$ and $\ip_\fk$ is non-degenerate it suffices to prove $\langle \omega(A,v,w),B\rangle=0$ for all $v,w\in V$ and $A,B\in \fk$. This equation follows from
\begin{eqnarray*}
\langle \omega(A,v,w),B\rangle&=& \langle [[A,v],w],B\rangle+\langle[[v,w],A],B\rangle+\langle[[w,A],v],B\rangle\\
&=&\langle [B,[A, v]],w\rangle+\langle [[A,B],v],w\rangle -\langle [A,[B,v]],w\rangle \ =\ 0.
\end{eqnarray*}
It remains to show $\omega|_{V\times V\times V}=0$. For this we need a few preparations. 
By (\ref{Tuu}) we have $T^{-1}[B,T_u]=[B,T_u](0)=dB(T_u(0))-dT_u(B(0))$ for $B\in\fk_+$ and $u\in V_-$. Since $B\in\fk_+$ implies $B(0)=0$ this yields $T^{-1}[B,T_u]=dB(u)=[B,u]$. Applying $T$ we obtain 
\begin{equation}\label{na}
[B,T_u]=T_{[B,u]}.
\end{equation}
Using this equation and the fact that $T$ is an isometry we see that
$$\langle [T_u,T_v],B\rangle=\langle T_u,[T_v,B]\rangle=-\langle T_u,T_{[B,v]}\rangle=-\langle u,[B,v]\rangle=\langle [u,v],B\rangle$$ 
holds for all $u,v\in V_-$ and $B\in\fk_+$. Since $\fk_+$ is non-degenerate we get
\begin{equation}\label{und}
[T_u,T_v]=[u,v]
\end{equation}
for all $u,v\in V_-$.
Now take $v\in V_-,\eta\in V_+$. Then 
$$\langle [v,\eta], T_w\rangle=-\langle v,[ T_w,\eta]\rangle=-\langle v, dT_w\eta\rangle=\langle v,A_\eta w\rangle =\langle A_\eta v,w\rangle =\langle T_{A_\eta v},T_w\rangle$$
holds for all $w\in V_-$. Since $V_-$ is non-degenerate we obtain
\begin{equation} \label{is}
[v,\eta]=T_{A_\eta v}\in\fk_-.
\end{equation}

Now we proceed with the proof of the Jacobi identity.
Take $u,v,w\in V_-$. By (\ref{R}) and (\ref{und}) the condition $\omega(u,v,w)=0$ is equivalent to first Bianchi identity for $R^M$.  

Now take $u,w\in V_-$, $\eta\in V_+$. From (\ref{ET}), (\ref{R}), (\ref{und}) and (\ref{is}) we get
\begin{eqnarray*} 
\omega(u,w,\eta) &=&[[u,w],\eta]+[[w,\eta],u]+[[\eta,u],w]\\
&=&-R^\perp(u,w)(\eta)+\alpha(A_\eta w,u)-\alpha(A_\eta u,w),
\end{eqnarray*}
which vanishes by the normal Gauss equation.

For $w\in V_-$, $\zeta,\eta\in V_+$ we have $\omega(w,\zeta,\eta)\in V_-$. Hence $\omega(w,\zeta,\eta)=0$ holds if and only if $\langle \omega(w,\zeta,\eta), v\rangle =0$ for all $v\in V_-$. Since $\langle \omega(w,\zeta,\eta), v\rangle=-\langle\eta,\omega(v,w,\zeta) \rangle =0$ by the above considerations we get  $\omega(w,\zeta,\eta)=0$. 

Instead of using (\ref{R}) Jacobi identity on $V\times V\times V$ except of $V_+\times V_+\times V_+$ can be proven only by (\ref{ET}) as done in \cite{EH}.

Finally we have to check that $\omega=0$ holds on $V_+\times V_+\times V_+$. Since $M$ is full we have $[\fk_-,V_-]= V_+$. Hence it suffices to prove $\omega(\zeta,\eta, [A,v])=0$ for all $\zeta,\eta\in V_+$, $A\in\fk_-$ and $v\in V_-$. Actually we need only ${\rm codim} _{V_+} [\fk_-,V_-]\le 2$ at this point. Because of
\begin{eqnarray*}
[[\zeta,\eta],[A,v]]&=&[[[\zeta,\eta],A],v]+[A,[[\zeta,\eta],v]]\\
&=&[[[\zeta,A],\eta]+[\zeta,[\eta,A]],v]+[A,[[\zeta,v],\eta]+[\zeta,[\eta,v]]]\\[1ex]
[[\eta,[A,v]],\zeta]&=&[[[\eta,A],v],\zeta]+[[A,[\eta,v]],\zeta]\\[1ex]
[[[A,v],\zeta],\eta]&=&[[[A,\zeta],v],\eta]+[[A,[v,\zeta]],\eta]
\end{eqnarray*}
we have
\begin{eqnarray*}
\omega(\zeta,\eta, [A,v])&=&\ \ [[\eta,[A,\zeta]],v]+[[[A,\zeta],v],\eta]+ [[[A,\eta],\zeta],v]+ [[v,[A,\eta]],\zeta]\\
&&+[A,[[\zeta,v],\eta]]+[\eta,[A,[\zeta,v]]]+[A,[\zeta,[\eta,v]]]+[\zeta,[[\eta,v],A]]\\[1ex]
&=&-[[v,\eta],[A,\zeta]]-[[\zeta,v],[A,\eta]]\\
&&+[A,[[\zeta,v],\eta]]+[\eta,[A,[\zeta,v]]]+[A,[\zeta,[\eta,v]]]+[\zeta,[[\eta,v],A]]\\[1ex]
&=&\omega(A,\zeta,[v,\eta])+\omega(A,\eta, [\zeta,v])\ =\ 0.
\end{eqnarray*}
This finishes the proof of the Jacobi identity.

Finally let us show that $D$ satisfies the conditions from Def.~\ref{Dq}. 
 It is easy to see that $D$ is antisymmetric. Indeed, for $D|_{\fk_+\oplus V_+}=0$ this is obvious and for $v\in V_-$, $A\in\fk_-$ we have
 $$\langle Dv,A\rangle = -\langle Dv,DDA\rangle=-\langle v, DA\rangle$$ since $D_{V_-}=T$ is an isometry.
 
$D$ is a derivation. Indeed, the restriction of $D$ to the subalgebra $\fk_+\oplus V_+$is a derivation since it is zero there.  By (\ref{ET}) we have
\begin{eqnarray*}
D[u,v]&=&0\ =\ \alpha(u,v)-\alpha(v,u)\ =\ dT_u(v)-dT_v(u)\ =\ [T_u,v]+[u,T_v]\\ &=&[Du,v]+[u,Dv]
\end{eqnarray*}
for $u,v\in V_-$.
Using (\ref{ET}) and (\ref{is}) we get 
\begin{eqnarray*}
D[u,\eta]&=&DT_{A_\eta u}\ =\ D^2(A_\eta u)\ =\ -A_\eta u\ =\ dT_u(\eta)\ =\ [T_u,\eta]\ =\ [Du,\eta]\\ &=&[Du,\eta]+[u,D\eta]
\end{eqnarray*}
and
$$D[T_u,\eta]=-D\,A_\eta u=-T_{A_\eta u}=-[u,\eta]=[DDu,\eta]=[DT_u,\eta]=[DT_u,\eta]+[T_u,D\eta]$$
for  $u\in V_-$ and $\eta\in V_+$. By (\ref{und}) we have
$$D[T_u,v]=0=-[u,v]+[T_u,T_v]=[DT_u,v]+[T_u,Dv]$$
for all $u,v\in V_-$. Because of (\ref{na}) we have
$$D[B,u]=T_{[B,u]}=[B,T_u]=[B,Du]=[DB,u]+[B,Du]$$
and
$$D[B,T_u]=DT_{[B,u]}=-[B,u]=[B,DT_u]=[DB,T_u]+[B,DT_u]$$
 for $B\in\fk_+$, $u\in V_-$.
Moreover, by (\ref{und})
\begin{eqnarray*}
D[T_u,T_v]&=&0\ =\ -D[u,v]\ =\ -[Du,v]-[u,Dv]\ =\ -[T_u,v]-[u,T_v]\\ &=&[T_u,DT_v]+[DT_u,T_v]
\end{eqnarray*}
for all $u,v\in V_-$.

 Let $\tau_D$ be the involution associated with $D$. 
Using the notation introduced in (\ref{not}) we get $\fg_+^+=\fk_+$ and 
$\fg_+^-=\fk_-$. This  gives $[\fg_+^-,\fg_+^-]=\fg_+^+$ 
since $\fk$ is the Lie algebra of the transvection group of the symmetric space $M$.  

The equality $M=M_{\fg,\Phi}$ holds by construction. 
\qed 

\begin{re}\label{fe}{\rm
{}From Remark~\ref{R4.2.} we can see that $M_{\fg,\Phi}$ is normal if $D$ is inner. Now we will show the converse. Assume that $M:=M_{\fg,\Phi}$ is normal. Then there exists an element $x_0\in \bigcap_{x\in M} \cN_xM.$ Identifying as usual the affine space $V$ with the vector space $V$ we may write $x_0=-\xi$. Since $x_0$ is a fixed point of all reflections we get 
$$0=T_u(-\xi) =dT_u(-\xi)+T_u(0)=-[T_u,\xi]+u=-[Du,\xi]+u$$
for all $u\in V_-$, hence $\ad(\xi)=-D^{-1}=D$ on $\fg_+^-$, which implies $D=\ad(\xi)$ on whole $\fg$. 

In particular this shows that in the case where the assumptions for Ferus' construction or for Naitoh's version are fulfilled these constructions yield the same object as our construction. Indeed, let $M\subset V$ be full and normal and let $(\fg,\ip,\Phi)$ be as constructed above. Choose $x_0$ and $\xi$ as just explained. Now take $x_0$ as the origin of $V$ and identify the affine with the vector space $V$ according to this choice as done by Ferus. Then $\xi$ is in $M$. Taking $\xi$ as a base point of $M$ Ferus' construction yields the same Lie algebra, the same Cartan decomposition and the same inner derivation as our.
}
\end{re}

\section{Weak triples as extensions of non-degenerate ones}
In the forthcoming paper \cite{Kext2} we will see that (non-degenerate) extrinsic symmetric triples can be described by a means of a certain cohomology set.  This description can be used to get explicit classification results in various concrete situations. Therefore we now want to clarify the relation between extrinsic symmetric triples and weak extrinsic symmetric triples. Let  $(\fg,\ip,\Phi)$, $\Phi=(D,\theta)$, be a weak extrinsic symmetric triple. Let us consider the metric radical  $R:=\fg\cap \fg^\perp\subset\fg_-^+$ of $\fg$. Then $R$ is an ideal of $\fg$. The inner product $\ip$ induces a non-degenerate inner product on $\fg/R$. Furthermore, $\Phi$ induces a corresponding structure on $\fg/R$. We will denote $\fg/R$ together with these induced structures just by $\fg/R$.
\begin{pr} 
\begin{enumerate}
\item $R$  is contained in the centre of $\fg$.
\item The quotient $\fg_0:=\fg/R$ is an extrinsic symmetric triple.
\end{enumerate}
\end{pr}
\proof Since on the one hand $[R,\fg_+^-\oplus\fg_-]\subset[\fg_-^+, \fg_+^-\oplus\fg_-]\subset \fg_-^-\oplus\fg_+$ and on the other hand $[R,\fg]\subset R\subset \fg_-^+$ we obtain $[R,\fg_+^-\oplus\fg_-]=0$. Because of $\fg_+^+=[\fg_+^-,\fg_+^-]$ we get also $[R,\fg_+^+]=0$ by Jacobi identity.
\qed

If we choose a section $s:\fg_0\rightarrow \fg$, then 
$$\omega \in \Hom(\textstyle{\bigwedge^2} \fg_0,R),\quad \omega(X,Y)=[s(X),s(Y)]-s([X,Y])$$
 defines a cocycle $\omega\in Z^2(\fg_0,R)$, where $\fg_0$ acts trivially on $R$. This cocycle satisfies 
\begin{equation}\label{omega}
\theta^*\omega=-\omega,\quad D \omega=0,
\end{equation}
where $D\omega(X,Y)=\omega(DX,Y)+\omega(X,DY)$. Conversely, given an extrinsic symmetric triple $(\fg_0,\ip_0,\Phi_0)$, $\Phi_0=(D_0,\theta_0)$, a real vector space $R$ (considered as a trivial $\fg_0$-module) and a cocycle $\omega\in Z^2(\fg_0,R)$ satisfying~(\ref{omega}) we can consider the extension $\fg$ of $\fg_0$ by $R$ defined by $\omega$. As a vector space $\fg$ equals  $\fg=\fg_0\oplus R$. If we put $\ip=\ip_0\oplus 0$ and $\Phi=(D_0\oplus 0, \theta_0\oplus-\Id)$ on $\fg_0\oplus R$ we obtain a weak extrinsic symmetric triple $(\fg,\ip,\Phi)$.

We can define actions of $\theta_0$ and $D_0$ on $H^2(\fg_0,R)$ by $\theta_0[\omega]= [\theta^*_0\omega]$ and $D_0[\omega]=[D_0\omega]$. We define
$$H^2(\fg_0,R)^{D_0}_-=\{a\in H^2(\fg_0,R)\mid \theta_0(a)=-a,\ D_0a=0\}$$ and obtain

\begin{pr}\label{Phol}
Isomorphism classes of weak extrinsic symmetric triples $(\fg,\ip,\Phi)$ with metric radical $R$  for which $\fg/R$ is isomorphic to the extrinsic symmetric triple $(\fg_0,\ip_0,\Phi_0)$ correspond to elements of $H^2(\fg_0,R)^{D_0}_-/(\Aut(\fg_0,\ip_0,\Phi_0)\times \GL(R))$.
\end{pr}

\begin{re}{\rm The extension $\fg$ of $\fg_0$ by $R$ is full  if and only if $\fg_0$ is full and $\omega(\fg_0,\fg_0)=R$.}
\end{re}

\begin{re}\label{REK}{\rm 
Eschenburg and Kim noticed in \cite{EK}  that for each extrinsic symmetric space $M$ in a degenerate inner product space $V$ the map $\pi:V\rightarrow\bar V:= V/(V\cap V^\perp)$ restricts to an extrinsic symmetric immersion $\pi|_M:M\hookrightarrow \bar V$. In order to reduce the classification problem of extrinsic symmetric spaces in degenerate inner product spaces to those in non-degenerate ones they posed the problem of determining all extrinsic symmetric spaces $N\hookrightarrow W$ for which $N\hookrightarrow \bar W:=W/(W\cap W^\perp)$ is (extrinsic) isometric to $M\hookrightarrow \bar V$. Now we can answer this question. We may restrict ourselves to the case where $M$ and $N$ are full. Then $M\hookrightarrow \bar V$ is also full. Hence it is associated with a full extrinsic symmetric triple $(\fg_0,\ip_0,\Phi_0)$. Now we see that the isometry classes of all possible $N\hookrightarrow W$ correspond exactly to the isomorphism classes of full weak extrinsic symmetric triples $(\fg,\ip,\Phi)$ for which $\fg/R$ is isomorphic to $(\fg_0,\ip_0,\Phi_0)$. These can be classified according to Prop.~\ref{Phol}}.
\end{re}
 
\section{Examples}                                        %
In this section we will give some examples of extrinsic symmetric triples and associated extrinsic symmetric spaces, which illustrate the facts proven above.  

\begin{ex}[Parabola]\label{ExPar}{\rm 
Let us start with a very simple example of an extrinsic symmetric space 
that is contained in a degenerate affine subspace but in no proper 
non-degenerate subspace. 

Let $\fg$ be the three-dimensional Heisenberg algebra $\fg$. Let $e_1,e_2,e_3$ be a basis of $\fg$ such that 
$$[e_1,e_2]=[e_2,e_3]=0,\quad [e_1, e_3]=e_2.$$ Let $\ip$ be the inner product for which $e_1,e_2,e_3$ are orthogonal to each other and $$\langle e_1,e_1\rangle=\langle e_3,e_3\rangle=1,\quad \langle e_2,e_2\rangle=0.$$
We define $\theta$ by
$\fg_+=\RR\cdot e_1$ and $\fg_-=\Span\{ e_2, e_3\}$. Moreover, let $D$ be defined by
$$D(e_1)=e_3,\quad D(e_2)=0,\quad D(e_3)=-e_1.$$
Then $(\fg,\ip,\Phi)$  with $\Phi=(D,\theta)$ is a weak extrinsic symmetric triple. The associated extrinsic symmetric space is the parabola
$$\RR^{0,1} \ni t\longmapsto 
\gamma(t):=\exp(t\phi(e_1))(0)=-\textstyle {\frac12} 
t^2 e_2-te_3\subset \RR^{0,1,1}\subset\RR^{1,2}.$$ 
Here $\fg_0=\fg/R=\Span\{e_1,e_3\}$ and the associated extrinsic symmetric space is the line $\RR\cdot e_3$.
}\end{ex} 

\begin{ex}[Three-dimensional flat space] \label{Ex*} {\rm \  
This example will show that the transvection group of a non-full extrinsic symmetric space $M\hookrightarrow V$ can be larger than the transvection group of the symmetric space $M$. We start with the construction of an extrinsic symmetric triple $(\fg,\ip,\Phi)$. Let $e_1,\dots,e_6$ be the standard basis of $\RR^6$ and let $\sigma^1,\dots,\sigma^6$ be the dual basis of $(\RR^6)^*$. Furthermore, let $b_1,b_2$ be a basis of $\RR^2$. Let $\omega\in Z^2(\RR^6,\RR^2)$ be the cocycle 
$$\omega=(\sigma_1\wedge\sigma_5+\sigma_2\wedge\sigma_4)\otimes b_1+(\sigma_1\wedge\sigma_6+\sigma_3\wedge\sigma_4)\otimes b_2.$$ Now let $\fg$ be the central extension of $\RR^6$ by $\RR^2$ defined by $\omega$. Let $\ip$ be the inner product on $\fg$ whose restriction to $\RR^6$ is the standard inner product and for which $b_1$ and $b_2$ span the metric radical. Let $\theta$ be defined by
$$\fg_+=\Span\{e_1,e_2,e_3\},\quad \fg_-=\Span\{e_4,e_5,e_6,b_1,b_2\}.$$
Moreover, we define $D$ by
$$D(b_1)=D(b_2)=0,\quad D(e_i)=-e_{i+3},\ i=1,2,3,\quad D(e_i)=e_{i-3},\ i=4,5,6.$$

We consider the vector space $V=\RR^7$ together with the scalar product of signature $(2,5)$ defined by
$$\langle x,y\rangle =x_1y_1+x_2y_2+x_3y_3+x_4y_6+x_5y_7+x_6y_4+x_7y_5$$
for $x,y\in\RR^7$. Then the map which sends $(x_1,\dots,x_5)\in\fg_-\cong \RR^5$ to $(x_1,\dots,x_5,0,0)\in V$ is an isometric embedding. 

Then an easy computation shows that the extrinsic symmetric space which is associated with $(\fg,\ip,\Phi)$ equals
$$f:M=\RR^3\hookrightarrow \fg_-\hookrightarrow V,\quad (t,s,r)\longmapsto (t,s,r,st,rt,0,0).$$
We have
\begin{equation}\label{ds}
ds_{f(t,s,r)}=
\mbox{\begin{small}$
\left(\begin{array}{cccccccc}
-1&0&0&0&0&-2s&-2r\\
0&-1&0&0&0&-2t&0\\
0&0&-1&0&0&0&-2t\\
-2s&-2t&0&1&0&-2(s^2+t^2)&-2rs\\
-2r&0&-2t&0&1&-2rs&-2(r^2+t^2)\\
0&0&0&0&0&1&0\\
0&0&0&0&0&0&1
\end{array}\right).
$\end{small}}
\end{equation}
{}From this one easily computes 
\begin{equation}\label{s}
s_{f(t,s,r)}(f(\bar t,\bar s,\bar r))=f(2t-\bar t, 2s -\bar s, 2r -\bar r),
\end{equation}
thus $s_x(M)=M$ for all $x\in M$. Now we use the reflections 
$$s_0:= s_{f(0,0,0)},\quad s_1:=s_{f(0,u,0)},\quad s_2:=s_{f(0,0,u)},\ u\in\RR$$ to define the transvection $\ph:=(s_1\circ s_0\circ s_2)^2$.  Then $\ph$ is the identity on $M$ by (\ref{s}). However,  (\ref{ds}) shows that  $d\ph$ is not the identity on $\RR^7$ if $u\not=0$. Hence the transvection group of the extrinsic symmetric space $M\subset\RR^7$ is larger than the transvection group of $M$.
}
\end{ex}
The following examples are constructed by a general principle called quadratic extension. Here we will use only a `light version' of this extension procedure. 
A detailed description of the general construction, which will allow a systematic study of extrinsic symmetric triples will follow in \cite{Kext2}. 
Let $(\fl,\theta_\fl)$ be a Lie algebra with involution, let $(\fa,\ip_\fa)$ be a pseudo-Euclidean vector space and $\theta_\fa$ an involutive isometry on $\fa$. Moreover, let $\rho$ be an orthogonal representation of $\fl$ on $\fa$ such that $\theta_\fa\circ\rho(\theta_\fl(L))=\rho(L)\circ \theta_\fa$. 
Let $\fg$ be the extension of the semi-direct sum $\fa\rtimes_\rho\fl$ by the abelian Lie algebra $\fl^*$ defined as follows. We identify $\fg$ as a vector space with $\fl^*\oplus\fa\oplus\fl$ and put
$$[L+A,Z]=\ad^*(L)(Z),\quad \proj_{\fl_*}[A_1+L_1,A_2+L_2]=\langle \rho(\cdot)A_1,A_2\rangle $$
for all $Z\in\fl^*$, $A,A_1,A_2\in\fa$, $L,L_1,L_2\in\fl$.  Let $\ip$ be the unique inner product on $\fg$ such that $\fl$ and $\fl^*$ are isotropic, $\fa\perp \fl^*\oplus\fl$,  $\ip|_{\fa\times \fa}=\ip_\fa$ and $\ip|_{\fl^*\times \fl}$ is the dual pairing of $\fl^*$ and $\fl$. Furthermore, we put
$\theta_\fg=\theta_\fl^*\oplus\theta_\fa\oplus\theta_\fl$. Then it is not hard to prove that $\dd:=(\fg,\ip,\theta)$ is a metric Lie algebra with involution, see \cite{KO2}.

\begin{ex}[Three-dimensional Cahen-Wallach spaces I]\label{Ex3}{\rm 
We construct examples of solvable Lorentzian extrinsic symmetric spaces. 

We consider $\fg=\dd$ with 
\begin{eqnarray*} 
&&\fl=\fsl(2,\RR)=\{[H,X]=2Y, [H,Y]=2X, [X,Y]=2H\}, \\
&&\mbox{with }\theta_\fl \mbox{ given by } \fl_+=\RR\cdot H,\ 
\fl_-=\Span\{X,Y\},\\ 
&&(\rho,\fa,\ip_\fa,\theta_\fa)=(\ad,\fl,\ip_\fl,-\theta_\fl), \mbox{ 
where $\ip_\fl$ is 
the Killing 
form on $\fl$}. 
\end{eqnarray*} 
For $\xi:=\frac 12 X$ we have $\ad(\xi)^3=-\ad(\xi)$ and   $(\fg,\ip,\theta)$ together with $D:=\ad (\xi)$ is an 
extrinsic symmetric triple. Let us study the associated extrinsic symmetric space in 
more detail. If 
$a_X,a_Y,a_H$ denote the images of $X,Y,H$ under the isomorphism 
$\fa\cong\fl$ and if $\sigma_X,\sigma_Y,\sigma_H$ is a basis of $\fl^*$ dual to $X,Y,H$, 
then we have 
$$\fg_+^+=\Span\{ a_X\},\ \fg_+^-=\Span\{\sigma_H, a_Y, H\}, \fg_-=\Span\{\sigma_X,\sigma_Y, a_H, X,Y\}$$ 
Hence $(\fg_{+},\ip|_{\fg_{+}\times\fg_{+}},\tau_D|_{\fg_{+}})$ is a 
solvable non-abelian symmetric triple. The associated symmetric 
space $M$ is 3-dimensional and has Lorentzian signature, thus it is one of 
the 3-dimensional Cahen-Wallach spaces. Here it is embedded in 
$\fg_{-}\cong\RR^{2,3}$. One can also obtain the other three-dimensional Cahen-Wallach space as an extrinsic symmetric space if one uses $\fl=\fsu(2)$ 
instead of $\fsl(2,\RR)$.

Let us compute the shape operator of $M\hookrightarrow \RR^{2,3}$. We 
denote the dual basis of $X,Y,H$ by $\sigma_{X},\sigma_{Y}, 
\sigma_{H}$. Using that the isometry $T^{-1}:\fg_{+}^{-}\rightarrow \fg_{-}^{-}$ is given 
by $T^{-1}=-\ad (\xi)|_{\fg_{+}^{-}}$ we see that the mean curvature vector $h$ at $\xi\in M$ equals 
\begin{eqnarray*}
h&=&\textstyle{\frac13}(2\alpha(\sigma_{Y},Y)+\textstyle{\frac18} 
\alpha(a_H, a_H)) \ 
=\ \textstyle{\frac13}(2[T_{\sigma_{Y}},Y]+\textstyle{\frac18} 
[T_{a_H}, a_H]) \\
&=& \textstyle{\frac13}(2[\sigma_{H},Y]+\textstyle{\frac18} 
[-a_Y, a_H]) =-2\sigma_{X}.
\end{eqnarray*}
Since $A_h(u)=-T_u(h)=2T_u(\sigma_{X})$ we obtain 
$$A_{h}(Y)=-4\sigma_{Y},\ A_{h}(\sigma_{Y})=A_{h}(a_{H})=0.$$ 
In particular, $A_{h}$ is 2-step nilpotent. 

The submanifold $M\hookrightarrow \RR^{2,3}$ is full since $[\fg_{+}^{-}, \fg_{-}^{-}]=\fg_{-}^{+}$. Moreover, it is normal since $D$ is inner. }\end{ex} 

\begin{ex}[Three-dimensional Cahen-Wallach spaces II]{\rm 
Here we will describe a three-dimensional Cahen-Wallach space as an extrinsic symmetric space of $\RR^{3,3}$ which is full but not normal. In particular it is not of Ferus type in the sense of \cite{EK}.
Now we take $\fg=\dd$ with 
\begin{eqnarray*} 
&&\fl=\fsl(2,\RR) \mbox{ with basis } H,X,Y, \\
&&\theta_\fl \mbox{ given by } \fl_+=\RR\cdot H,\ 
\fl_-=\Span\{X,Y\},\\ 
&&(\rho,\fa)=(\rho_2\oplus\rho_2,\fa_2\oplus\fa_2\cong \fa_2\otimes \RR^2), \mbox{ 
where $(\rho_2,\fa_2)$ is the standard representation}\\&& \mbox{of $\fsl(2,\RR)$. We identify $\fa_2\oplus\fa_2\cong \fa_2\otimes\RR^2$ and consider}\\
&&\ip_\fa =\mbox{\small{$\left(\begin{array}{cc}0&-1\\1&0\end{array}\right)
\otimes \left(\begin{array}{cc}0&-1\\1&0\end{array}\right),$}}\quad
\theta_\fa =\mbox{\small{$\left(\begin{array}{cc}1&0\\0&-1\end{array}\right)
\otimes \left(\begin{array}{cc}1&0\\0&-1\end{array}\right)$}}.
\end{eqnarray*} 
We define an anti-symmetric map $D_\fa:\fa\rightarrow\fa$ by 
$$ D_\fa(b_1)=D_\fa(b_2)=0,\ D_\fa(b_3)=b_4,\ D_\fa(b_4)=-b_3$$
and a derivation $D_\fl:\fl\rightarrow\fl$ by
$$ D_\fl(X)=0,\ D_\fl(Y)=H,\ D_\fl(H)=-Y.$$
Then $D:=-(D_\fl)^*\oplus D_\fa\oplus D_\fl$ is an antisymmetric derivation on $\dd$ satisfying  $D^3=-D$ and it is easy to check that $(\fg,\ip,\theta)$ together with $D$ is an extrinsic symmetric triple. Let $a_1,a_2$ and $e_1,e_2$ denote the standard bases of $\fa_2$ and $\RR^2$, respectively. Then
$$
\begin{array}{ll}
b_1:=(1/\sqrt 2)\cdot (a_1\otimes e_1 +a_2\otimes e_2), & b_2:= (1/\sqrt 2)\cdot(a_1\otimes e_2 -a_2\otimes e_1),\\[0.5ex]
b_3:=(1/\sqrt 2)\cdot (a_2\otimes e_2 -a_1\otimes e_1), & b_4:= (1/\sqrt 2)\cdot(a_1\otimes e_2 +a_2\otimes e_1)
\end{array}
$$
is a basis of $\fa$. We have
$$
\begin{array}{ll}
\fg_+^+=\Span\{b_1\}, & \fg_+^-=\Span\{\sigma_H,b_3,H \},\\[0.5ex]
\fg_-^+=\Span\{\sigma_X,b_2,X \} ,& \fg_-^-=\Span\{ \sigma_Y,b_4,Y\}.
\end{array}
$$
This yields $[\fg_+^-,\fg_-^-]=\Span \{\sigma_X,b_2,X\}$. Hence the associated extrinsic symmetric space $M\hookrightarrow\fg_-\cong \RR^{3,3}$ is full. It is not hard to see that $D$ is not an inner derivation. Thus $M$ is not normal in $\RR^{3,3}$. A direct calculation shows that $M$ can be parametrised by
\begin{eqnarray*}
 (r,s,t)&\longmapsto& \Big( \textstyle{ (\frac{2s^2}r -t) \sinh 2r +\frac{s^2}{r^2}(1-\cosh 2r),\  -(\frac{2s^2}r -t) \cosh 2r +\frac{s^2}{r^2}\sinh 2r,}\\
&& \textstyle{\frac sr(1-\cosh 2r),\  -\frac sr \sinh 2r,\  \frac12\cosh 2r-\frac12,\ \frac12\sinh 2r\Big).}
	\end{eqnarray*}
	According to \cite{EK} the shape operator $A_h$ with respect to the mean curvature vector must satisfy ${A_h}^2=0$. In fact, one computes $h=-\sigma_x$ and 
$$A_h(Y)= -2\sigma_Y,\ A_h(b_4)=A_h(\sigma_Y)=0$$	
at $0\in M$. 
} \end{ex}


\begin{thebibliography}{MMM} 
\bibitem[CP]{CP} M.\,Cahen, M.\,Parker, {\sl Pseudo-Riemannian symmetric spaces.} Mem.\,Amer. Math.\,Soc.\,{\bf 24} (1980), no.~229. 
\bibitem[EH]{EH} J.-H.\,Eschenburg, E.\,Heintze, {\sl Extrinsic symmetric 
spaces and orbits of s-representations.} Manuscripta Math.\,{\bf 88} (1995), 
517 -- 524. 
\bibitem[EK]{EK} J.-H.\,Eschenburg,  J.\,R.\,Kim, {\sl Indefinite extrinsic symmetric spaces}, preprint.
\bibitem[F1]{F1} D.\,Ferus, {\sl Produkt--Zerlegung von Immersionen mit 
paralleler zweiter Fundamentalform.} Math.\,Ann.\,{\bf 211} (1974), 1 -- 5. 
\bibitem[F2]{F2} D.\,Ferus, {\sl Immersions with Parallel Second Fundamental Form.} Math.\,Z.\,{\bf 140} (1974), 87-92. 
\bibitem[F3]{F3} D.\,Ferus, {\sl Symmetric Submanifolds of Euclidean Space.} Math.\,Ann.\,{\bf 247} (1980), 81 -- 93. 
\bibitem[Ka]{Kext2}I.\,Kath, {\sl Indefinite extrinsic symmetric spaces II}
\bibitem[KO1]{KO1} I.\,Kath, M.\,Olbrich, 
{\sl Metric Lie algebras and quadratic extensions.} Transform.\,Groups {\bf 11} (2006), no.~1, 87 -- 131. 
\bibitem[KO2]{KO2} I.\,Kath, M.\,Olbrich, 
{\sl On the structure of pseudo-Riemannian symmetric spaces.} 
arXiv:math.DG/0408249, 2004. 
\bibitem[KO3]{KOesi} I.\,Kath, M.\,Olbrich,  {\sl The classification problem for 
pseudo-Riemannian symmetric spaces} in Recent developments in 
pseudo-Riemannian Geometry,
ESI Lectures in Mathematics and Physics, EMS Publishing House, 2008, 1 -- 52. 
\bibitem[Ki]{Kim} J.\,R.\,Kim, {\sl Indefinite Extrinsic Symmetric Spaces.}
Dissertation, Augsburg 2005, Shaker Verlag. 
\bibitem[Ko]{Ko} S.\,Kobayashi, {\sl Isometric imbeddings of compact symmetric spaces.} T\^ohoku Math.\,J. {\bf 20} (1968), 21 -- 25.
\bibitem[KoNa] {KoNa} S.\,Kobayashi, T.\,Nagano, {\sl On filtered Lie algebras and geometric structures.} I.\,J.\,Math. Mech. {\bf 13} (1964), 875 -- 907.
\bibitem[Na]{Na} T.\,Nagano, {\sl Transformation groups on compact symmetric spaces.} Trans. Amer. Math. Soc. {\bf 118} (1965), 428 -- 453. 
\bibitem[Nai]{Nai} Naitoh,\,H., {\sl Pseudo-Riemannian symmetric $R$-spaces.} Osaka J.\,Math. {\bf 21} (1984), 733 -- 764.
\bibitem[S]{S} W.\,Str\"ubing, {\sl Symmetric submanifolds of Riemannian 
manifolds.} Math.\,Ann. {\bf 245} (1979), 37 -- 44. 
\end{thebibliography}
\end{document}